\documentclass[11pt]{article}

\usepackage{epsfig,amsmath,latexsym,amssymb}
\usepackage{mathtools}
\usepackage{graphicx}
\usepackage{lscape}
\usepackage{picture, eso-pic, tikz} 
\usepackage{dsfont}
\usepackage{listings}
\usepackage{changes}
\usepackage{hyperref}
\usepackage{bbding}

\renewcommand{\baselinestretch}{1.1}
\def\R{{\mathbb R}}  
\def\p{{\mathbb P}}  
\def\E{{\mathbb E}}  %
\def\Beweis{\footnotesize}

\newcommand{\Remm}[1]{}
\newtheorem{theo}{Theorem}[section]

\newtheorem{prop}[theo]{Proposition}

\newtheorem{cor}[theo]{Corollary}
\newtheorem{defi}[theo]{Definition}
\newtheorem{model ass}[theo]{Model Assumptions}

\def\EndProof{\hfill {\scriptsize $\Box$}}

\numberwithin{equation}{section}

\definecolor{MyGray}{rgb}{0.92,0.92,0.92}


\newcommand{\bl}[1]{\textcolor{blue}{{#1}}}

\definecolor{British racing}{rgb}{0.0, 0.5, 0.0}

\def\bS{\boldsymbol{S}}
\def\bT{\boldsymbol{T}}

\def\bX{\boldsymbol{X}}

\def\b0{\boldsymbol{0}}

\def\br{\boldsymbol{r}}

\def\b0{\boldsymbol{0}}

\usepackage{todonotes}
\newcommand{\Comments}{1}
\newcommand{\mynote}[2]{\ifnum\Comments=1\textcolor{#1}{#2}\fi}
\newcommand{\mytodo}[2]{\ifnum\Comments=1%
  \todo[linecolor=#1!80!black,backgroundcolor=#1,bordercolor=#1!80!black]{#2}\fi}

\ifnum\Comments=1               
\fi

\ifnum\Comments=1               
\fi

\ifnum\Comments=1               
\fi

\begin{document}
\author{Mario V.~W\"uthrich\footnote{RiskLab, Department of Mathematics, ETH Zurich, mario.wuethrich@math.ethz.ch}
}

\date{Version of \today}
\title{\vspace{-.5cm}Auto-Calibration Tests for Discrete Finite Regression Functions}
\maketitle

\vspace{-.4cm}
\begin{abstract}
\noindent  
Auto-calibration is an important property 
of regression functions for
actuarial applications. Comparably little is known about statistical testing of auto-calibration. Denuit et al.~(2024) recently published a test
with an asymptotic distribution that is not fully explicit and its evaluation needs non-parametric Monte Carlo sampling. In a simpler set-up, we present three test statistics with fully known and interpretable asymptotic distributions.

\medskip

\noindent
{\bf Keywords.} Auto-calibration, concentration curve, Lorenz curve, area between the curves.
\end{abstract}

\section{Introduction}
\label{Section Introduction}
Recent actuarial and financial literature acknowledges the importance of the statistical
concept of auto-calibration; see, e.g., Kr\"uger--Ziegel \cite{Ziegel},  Denuit et al.~\cite{DenuitEtAl} and W\"uthrich \cite{W_Gini}. Select an integrable 
response variable $Y$ and
covariates $\bX$ with support ${\cal X}$.
\begin{defi}
A measurable regression function $\pi:{\cal X}\to \R$ is auto-calibrated for $(Y,\bX)$ if
\begin{equation*}
\pi(\bX) = \E \left[ \left. Y \right| \pi(\bX) \right],
\qquad \text{$\p$-a.s.}
\end{equation*}
\end{defi}
In an actuarial pricing context this means that every price cohort $\pi(\bX)$ is on average self-financing for the claims $Y$, or in other words, there is no systematic cross-financing within a pricing scheme designed by the regression function $\pi$.

Surprisingly, there is no mature literature on testing for auto-calibration. Most proposals only consider binary responses, e.g., 
Gneiting--Resin \cite{GneitingResin} discuss a bootstrap test 
and Dimitriadis et al.~\cite{Dimitriadis} study calibration bands. Recently, Denuit et al.~\cite{Denuit} presented an auto-calibration test that studies the difference between the concentration curve (CC) and the Lorenz curve (LC). 
Also this test requires simulations because the asymptotic 
distribution of the test statistics is not sufficiently explicitly. We take one step back here, and we present simpler test statistics with fully known
and interpretable asymptotic distributions, though, in a simpler set-up.

One needs three ingredients for an auto-calibration test. (a) A regression function $\pi: {\cal X} \to \R$. This regression function $\pi$ can be fully general, i.e., we do not require that it is close (in some metric) to the conditional mean $\E[Y|\bX]$, nor do we specify whether $\pi$ has been estimated from past data ${\cal D}$ or whether it has been set by an expert. (b) A pair $(Y,\bX)$. For simplicity, we assume that the response $Y$ is positive and square integrable. The covariates $\bX$ have support ${\cal X}$. (c) An i.i.d.~sample ${\cal T}=(Y_i, \bX_i)_{i\ge 1}$ for testing.
This sample should have the same law as $(Y, \bX)$. These three ingredients (a)-(c) are sufficient for testing for auto-calibration of $\pi$ for $(Y,\bX)$; if $\pi$ has been estimated from past data ${\cal D}$, we generally assume that $(Y,\bX)$, ${\cal T}$ and ${\cal D}$
are independent, and all subsequent statements need then be understood conditional on ${\cal D}$.

\section{Tests for auto-calibration}
Assume that the regression function $\pi:{\cal X} \to \R$ takes only finitely many (ordered) values $-\infty < \pi_1 < \cdots < \pi_K <\infty$.
This gives us a partition of the covariate space ${\cal X}$ with
\begin{equation}\label{finite partition case}
\p \left[ \pi(\bX) =\pi_k\right]=p_k>0
\qquad \text{ for all $1\le k \le K$.}
\end{equation}
We assume probabilities $p_k>0$, otherwise the corresponding part of the covariate space can be dropped.
In this finite partition case \eqref{finite partition case}, auto-calibration of $\pi$ for $(Y,\bX)$ is equivalent to
\begin{equation*}
\pi_k = \E \left[ \left. Y \right| \pi(\bX)=\pi_k \right]
\qquad \text{ for all $1\le k \le K$.}
\end{equation*}
Using the tower property, 
auto-calibration of $\pi$ for $(Y,\bX)$ implies for all $1\le k \le K$
\begin{equation}
S^{(k)}:=\E \left[ \left(Y-\pi(\bX) \right) \mathds{1}_{\{\pi(\bX)=\pi_k \}}\right]
=\label{auto-calibration test 1}
\E \left[ \left(\E \left[\left.Y\right|\pi(\bX)\right]-\pi(\bX) \right) \mathds{1}_{\{\pi(\bX)=\pi_k \}}\right]=0,
\end{equation}
this statement is essentially the same as W\"uthrich \cite[Proposition 4.1]{W_Gini}. For a given i.i.d.~sample ${\cal T}=(Y_i, \bX_i)_{i=1}^n$,
this motivates the statistics
\begin{equation*}
S^{(k)}_n=\frac{1}{n} \sum_{i=1}^n \left(Y_i-\pi(\bX_i) \right) \mathds{1}_{\{\pi(\bX_i)=\pi_k \}} \qquad \text{ for $1\le k \le K$}.
\end{equation*}
Under auto-calibration of $\pi$ for $(Y,\bX)$, these empirical quantities $S^{(k)}_n$, $1 \le k \le K$, converge to zero, $\p$-a.s., as $n \to \infty$,
and we have the following central limit theorem.
\begin{prop} \label{lemma 1}
Under auto-calibration
of $\pi$ for $(Y,\bX)$
\begin{equation*}
\sqrt{n}
\left(
S^{(1)}_n, \ldots, S^{(K)}_n\right)^\top
~ \Longrightarrow \quad {\cal N}\left(0,\, {\rm diag}\left(p_1 \tau_1^2, 
\ldots, p_K \tau_K^2\right)\right)
\qquad \text{ as $n\to \infty$,}
\end{equation*}
with conditional variances $\tau_k^2 = {\rm Var}\left(\left.Y \right|\pi(\bX)=\pi_k \right)$ for
$1\le k \le K$.
\end{prop}
The proof of this proposition is standard and based on characteristic
functions.

\medskip

{\bf Test 1.}
Under the null hypothesis of $\pi$ being auto-calibrated for $(Y,\bX)$, \eqref{auto-calibration test 1} is a necessary condition for all $1\le k \le K$.
We test this against the alternative that there exists a 
$1 \le k \le K$ with $S^{(k)}\neq 0$. 
Under the null hypothesis, Proposition \ref{lemma 1} gives us for $s>0$ and $n$ large
\begin{equation}\label{test 1}
\p \left[ \max_{1\le k \le K}\sqrt{n} |S^{(k)}_n| \le s \right] = \p \left[ \bigcap_{1\le k \le K} \{| \sqrt{n}S^{(k)}_n| \le s\}
\right]  ~\approx~ \prod_{k=1}^K  
\left(2 \Phi\left(\frac{s}{\sqrt{p_k}\tau_k}\right)-1\right).
\end{equation}
Often, it is beneficial to test for the maximum of the normalized quantities
$\sqrt{n}|S^{(k)}_n|/(\sqrt{p_k}\tau_k)$, to have all terms on the same scale. This provides asymptotic limit
$(2 \Phi(s)-1)^K$.

\medskip

Denuit et al.~\cite[formula (2.4)]{Denuit} consider an aggregated version
of $S^{(k)}$. Namely, auto-calibration of $\pi$ for $(Y,\bX)$ implies for all $1\le k \le K$
\begin{equation}\label{auto-calibration test 2}
T^{(k)}:=\E \left[ \left(Y-\pi(\bX) \right) \mathds{1}_{\{\pi(\bX)\le \pi_k \}}\right]=0.
\end{equation}
For a given i.i.d.~sample ${\cal T}=(Y_i, \bX_i)_{i=1}^n$, this motivates the statistics
\begin{equation*}
T^{(k)}_n=\frac{1}{n} \sum_{i=1}^n \left(Y_i-\pi(\bX_i) \right) \mathds{1}_{\{\pi(\bX_i) \le \pi_k \}}=
 \sum_{j=1}^k S^{(j)}_n.
\end{equation*}
The following corollary is an immediate consequence of 
Proposition \ref{lemma 1}.
\begin{cor}\label{corollary 2} 
Under auto-calibration
of $\pi$ for $(Y,\bX)$ 
\begin{equation*}
\sqrt{n}
\left(
T^{(1)}_n, \ldots, T^{(K)}_n\right)^\top ~
 \Longrightarrow \quad {\cal N}\left(0, \left(\sum\nolimits_{j=1}^{\min\{k,m\}}
 p_j\tau^2_j \right)_{1\le k,m\le K}\right)
\qquad \text{ as $n\to \infty$.}
\end{equation*}
\end{cor}
Thus, the aggregated statistics $(T^{(k)}_n)_{1=k}^K$ can asymptotically be described by a random walk
\begin{equation}\label{random walk}
Z_k = \sum_{j=1}^k \sqrt{p_j}\,\tau_j\, \varepsilon_j,
\end{equation}
with i.i.d.~standard Gaussian innovations $\varepsilon_j \sim {\cal N}(0,1)$ for $1\le j \le K$. 

\medskip

{\bf Test 2.}
Under the null hypothesis of $\pi$ being auto-calibrated for $(Y,\bX)$, \eqref{auto-calibration test 2} is a necessary condition for all $1\le k \le K$.
We test this against the alternative that there exists a 
$1 \le k \le K$ with $T^{(k)}\neq 0$. 
Under the null hypothesis, Corollary \ref{corollary 2} gives us for $s>0$ and $n$ large
\begin{equation}\label{test 2}
\p \left[ \max_{1\le k \le K}\sqrt{n} |T^{(k)}_n| \le s \right] ~\approx~ 
\p \left[ \max_{1\le k \le K} | Z_k | \le s \right].
\end{equation}

\medskip

Up to one point discussed below, asymptotic approximation \eqref{test 2} gives an explicit explanation to the intractable limit in Denuit et al.~\cite[Proposition 3.1]{Denuit}. Namely, the asymptotic distribution
of the test statistics in \eqref{test 2} corresponds to the maximum of the random walk \eqref{random walk} 
whose increments are fully determined by the probabilities
$(p_k)_{k=1}^K$, given in \eqref{finite partition case}, and the
conditional variances $(\tau_k^2)_{k=1}^K$, given in Proposition \ref{lemma 1}. These two parameter sets can be determined from 
past data ${\cal D}$, being independent of the i.i.d.~sample ${\cal T}$, 
see discussion in Section \ref{Section Introduction}. The rejection area
is then received by (easy) random walk simulations involving only these
two (estimated) parameter sets $(p_k)_{k=1}^K$ and $(\tau_k^2)_{k=1}^K$.
This seems simpler than the non-parametric Monte Carlo
method used in Denuit et al.~\cite[Section 3.1]{Denuit}.

\section{Testing for the area between the curves}
The consideration of $T^{(k)}$
is motivated by the difference of the CC and the LC. 
Denote by $F^{-1}_{\pi(\bX)}$ the left-continuous generalized inverse of the distribution function $F_{\pi(\bX)}$ of $\pi(\bX)$. The difference between the CC
and the LC at probability level $\alpha \in [0,1]$ is defined by
\begin{equation*}
U(\alpha) =\E \left[ \left(\frac{Y}{\E[Y]}-\frac{\pi(\bX)}{\E[\pi(\bX)]}\right) \mathds{1}_{\{\pi (\bX)\le F^{-1}_{\pi(\bX)}(\alpha) \}}\right].
\end{equation*}
For a regression function $\pi$ with discrete finite range \eqref{finite partition case},
$U(\cdot)$ only takes $K+1$ different values in the 
cumulative probabilities $\alpha_k := \sum_{j=1}^k p_j$,
and we set $\alpha_0=0$. Namely, we
have
\begin{equation}\label{auto-calibration test 3}
U^{(k)}:=U(\alpha_k)=\E \left[ \left(\frac{Y}{\E[Y]}-\frac{\pi(\bX)}{\E[\pi(\bX)]}\right) \mathds{1}_{\{\pi(\bX)\le \pi_k \}}\right]
\qquad \text{ for $1\le k \le K$.}
\end{equation}
Under unbiasedness  $\E[\pi(\bX)]=\E[Y]$, we have
\begin{equation*}
U^{(k)} = \frac{1}{\E[Y]}\,T^{(k)}= \frac{1}{\E[\pi(\bX)]}\,T^{(k)}.
\end{equation*}
These normalized differences $U^{(k)}$ motivate the study
of $T^{(k)}$ under auto-calibration of $\pi$ for $(Y,\bX)$, which
implies the above unbiasedness. Denuit et al.~\cite[Proposition 3.1]{Denuit} do not exploit an auto-calibration test for $T^{(k)}$, but rather for $U(\alpha)$. Unfortunately, the normalized quantities $U(\alpha)$ and $U^{(k)}$
are more involved. 
For a given i.i.d.~sample ${\cal T}=(Y_i, \bX_i)_{i=1}^n$, 
consider
\begin{equation*}
U^{(k)}_{n}=\frac{1}{n} \sum_{i=1}^n \left(\frac{Y_i}{\overline{y}}-\frac{\pi(\bX_i)}{\overline{\pi}} \right) \mathds{1}_{\{\pi(\bX_i) \le \pi_k \}},
\end{equation*}
with $\overline{y}$ being the empirical mean of  $(Y_i)_{i=1}^n$ and
$\overline{\pi}$ the empirical mean of $(\pi(\bX_i))_{i=1}^n$. Dealing with 
$U^{(k)}_{n}$ instead of $T^{(k)}_{n}$ is more cumbersome because
of these normalizations. These normalizations are mainly motivated
by the fact that they imply that both the CC and the LC are calibrated
to 1 for $\alpha \uparrow 1$. In statistical modeling, this then allows one to perform model selection by selecting the model that has the most convex CC,
as a higher convexity implies better discrimination; see
W\"uthrich \cite{W_Gini}. Similarly, in economics,
a more convex LC indicates higher inequality in wealth distribution. 
However, for testing of auto-calibration this normalization seems not
justified, and we give preference to the simpler
unscaled quantity $T^{(k)}_{n}$.
Note that
\begin{eqnarray}\nonumber
\sqrt{n}
U^{(k)}_{n}&=&
\sqrt{n}\,\frac{1}{\overline{y}}\,
T^{(k)}_{n} +\sqrt{n}
\left(\frac{1}{\overline{y}}-\frac{1}{\overline{\pi}} \right)
\frac{1}{n} \sum_{i=1}^n \pi(\bX_i) \mathds{1}_{\{\pi(\bX_i) \le \pi_k \}}
\\&=&\label{Denuit convergence}
\sqrt{n}\,\frac{1}{\overline{y}}\,
T^{(k)}_{n} +\sqrt{n}\,\frac{1}{\overline{y}\,\overline{\pi}}
\left(\overline{\pi}-\overline{y} \right)
\frac{1}{n} \sum_{i=1}^n \pi(\bX_i) \mathds{1}_{\{\pi(\bX_i) \le \pi_k \}}.
\end{eqnarray}
Corollary \ref{corollary 2} and Slutsky's theorem
give weak convergence of the first term in \eqref{Denuit convergence} to 
$Z_k/\E[Y]$. For the second term in \eqref{Denuit convergence}, one
establishes weak convergence of $\sqrt{n}(\overline{\pi}-\overline{y})$, and the other terms are treated by Slutsky's theorem. Finally, one needs to
compute the covariance between the two
terms in \eqref{Denuit convergence} to get the asymptotic variance of 
$\sqrt{n}U^{(k)}_{n}$. This is doable, but cumbersome. Therefore,
we prefer to study the non-normalized quantities $T^{(k)}_{n}$.

\medskip

Based on $U(\alpha)$, Denuit et al.~\cite[formula (4.4)]{DenuitABC} introduced the area between the curves (ABC) as a model selection criterion. The ABC is
defined by 
\begin{equation*}
{\rm ABC} =
\int_{0}^1 U(\alpha)\,d\alpha= \int_{0}^1 \E \left[ \left(\frac{Y}{\E[Y]}-\frac{ \pi(\bX)}{\E[\pi(\bX)]}\right) \mathds{1}_{\{\pi(\bX)\le F^{-1}_{\pi(\bX)}(\alpha)\}}\right]  d\alpha.
\end{equation*}
Again, we prefer the unscaled version. Under the
discrete finite regression function, we have
\begin{equation*}
{\rm ABC}^\circ :=
\int_{0}^1 \E \left[ \left(Y-\pi(\bX)\right) \mathds{1}_{\{\pi(\bX)\le F^{-1}_{\pi(\bX)}(\alpha)\}}\right]  d\alpha
=
\sum_{k=1}^{K-1} p_{k+1}\, T^{(k)}.
\end{equation*}
For a given i.i.d.~sample ${\cal T}=(Y_i, \bX_i)_{i=1}^n$, 
this motivates the an integrated random walk statistics
\begin{equation*}
\widehat{\rm ABC}_n^\circ =
\sum_{k=1}^{K-1} p_{k+1} T_n^{(k)}
=
\sum_{k=1}^{K-1} p_{k+1} \sum_{j=1}^k S_n^{(j)}
=
\sum_{k=1}^{K-1} \left(1-\alpha_{k}\right)S_n^{(k)}.
\end{equation*}
Under auto-calibration of $\pi$ for $(Y,\bX)$, statistics
$\widehat{\rm ABC}_n^\circ$ converges to zero, $\p$-a.s.
Slightly modifying the terms, we propose the
following weighted $L^2$-norm statistics of the increments
\begin{equation}\label{test 3}
V^2_n
:= \sum_{k=1}^{K} \left(1-\alpha_{k-1}\right)(S_n^{(k)})^2,
\end{equation}
thus, the random walk increments $S_n^{(k)}$ with different signs cannot compensate
each other. 
\begin{cor}\label{corollary 3} 
Under auto-calibration
of $\pi$ for $(Y,\bX)$
\begin{equation*}
n\,V^2_n \quad
 \Longrightarrow \quad \sum_{k=1}^{K} \left(1-\alpha_{k-1}\right)p_k\,\tau_k^2\,\chi_k^2
\qquad \text{ as $n\to \infty$},
\end{equation*}
where $\chi_k^2$ are i.i.d.~$\chi^2$-distributed random variables 
with one degree of freedom.
\end{cor}

{\bf Test 3.} Under the above assumptions, we can test for auto-calibration of $\pi$ for $(Y,\bX)$ by exploiting the limiting distribution of Corollary \ref{corollary 3} numerically. As in Test 2, this limiting distribution
only depends on the two parameter sets $(p_k)_{k=1}^K$ and $(\tau_k^2)_{k=1}^K$. 

Dropping the weighting $1-\alpha_{k-1}$ in \eqref{test 3} and scaling the
individual terms $(S_n^{(k)})^2$ by $p_k \tau_k^2$ gives a $\chi^2$-test with $K$ degrees of freedom.

\section{Conclusions}
This letter considers statistical testing of auto-calibration. In the simplified
set-up of a discrete finite regression function, we provide
three different test statistics that have fully known asymptotic distributions under auto-calibration, see \eqref{test 1}, \eqref{test 2}
and Corollary \ref{corollary 3}. These three test statistics consider random walk increments, a random walk and an integrated random walk.
The three test statistics can be used for statistical testing of auto-calibration in our simpler set-up; Test 2 is a modified version of Denuit et al.~\cite[Proposition 3.1]{Denuit}. 

In this letter, we did not cover a study of the powers of these tests. This will depend on the kind of violation of auto-calibration; in fact, we believe that it is beneficial to normalize all random walk increments to unit variance in any of the three presented tests. Another open problem is to generalize these tests to arbitrary regression functions, this seems feasible for Tests 2 and 3.

\bigskip

{\small 
\renewcommand{\baselinestretch}{.5}

}

\newpage

\appendix
\setcounter{section}{19}
\section*{Supplementary}

\bigskip

\subsection*{Proofs}

\bigskip

{\Beweis 
{\bf Proof of Proposition \ref{lemma 1}.} 
Set $\bS_n=(S^{(1)}_n, \ldots, S^{(K)}_n)^\top$.
For $\br \in \R^K$, consider the characteristic function
\begin{eqnarray*}
\E \left[ \exp \left\{i \sqrt{n}\br^\top \bS_n \right\}\right]
&=&
\E \left[ \exp \left\{ \frac{i}{\sqrt{n}}\sum_{j=1}^n \sum_{k=1}^K r_k \left(Y_j-\pi(\bX_j) \right) \mathds{1}_{\{\pi(\bX_j)=\pi_k \}}\right\}\right]
\\&=&\prod_{j=1}^n\sum_{k=1}^K
\E \left[ \exp \left\{ \frac{i}{\sqrt{n}}  r_k \left(Y_j-\pi(\bX_j) \right) \right\}\mathds{1}_{\{\pi(\bX_j)=\pi_k \}}\right]
\\&=&
\exp \left\{ n \log
\left(\sum_{k=1}^K p_k\,
\E \left[\left. \exp \left\{ \frac{i}{\sqrt{n}}  r_k \left(Y-\pi(\bX) \right) \right\}\right| \pi(\bX)=\pi_k\right] \right) \right\}
\\&=&
\exp \left\{ n \log
\left(\sum_{k=1}^K p_k\left(1 -
\frac{r^2_k}{2n}\E \left[\left.  \left(Y-\pi(\bX) \right)^2
\right| \pi(\bX)=\pi_k\right] + o(n^{-1}) \right)\right) \right\}
\\&=& \prod_{k=1}^K
\exp \left\{-r_k^2\, \frac{p_k \tau_k^2}{2}   \right\}
\exp\{o(1)\} \qquad \text{ as $n\to \infty$,}
\end{eqnarray*}
where in the second last step we use auto-calibration
of $\pi$ for $(Y,\bX)$. This completes the proof.
\EndProof}

\bigskip

\newpage

\subsection*{Example}
We study a gamma distribution example with $K=6$ expected
response levels $(\pi_k)_{k=1}^K$. 
Table \ref{parameter values} shows the selected parameters. 
Firstly, we choose the probabilities $(p_k)_{k=1}^K$ such that the
boundary levels $\pi_k \in \{10,15\}$ receive the smallest probabilities,
and the levels in the middle $\pi_k \in \{12,13\}$ get the highest
probabilities. This is a quite
common feature in real data. Secondly, the variance parameters
$(\tau^2_k)_{k=1}^K$ are increasing in regression means 
$(\pi_k)_{k=1}^K$. Also this is a rather common feature, e.g., 
a Poisson or a gamma generalized linear model (GLM) have this property.
Based
on these parameters, we simulate first the regression level $\pi_k$ using
the probabilities $(p_k)_{k=1}^K$. Based on this level $\pi_k$, we
then simulate the response $Y|_{\pi(\bX)=\pi_k} \sim \Gamma(\gamma_k,c_k)$
with shape parameter $\gamma_k= 3 \pi_k$ and scale parameter $c_k=3$.
This gives us conditional mean $\pi_k$ and conditional variance $\tau_k^2=\pi_k/3$, see Table \ref{parameter values}. In particular, 
auto-calibration is fulfilled in this example because we simulate from the correct means.

\begin{table}[htb!]
\centering
{\small
\begin{center}
\begin{tabular}{|l||cccccc|}
\hline
$k$& 1 & 2 & 3 & 4 & 5 & 6\\\hline
\hline
$\pi_k$& 10 & 11 & 12 & 13 & 14 & 15\\
$p_k$ & $10/100$& $15/100$& $25/100$& $25/100$& $15/100$& $10/100$\\
$\tau^2_k$& 10/3 & 11/3 & 12/3 & 13/3 & 14/3 & 15/3\\
                                                         \hline
\end{tabular}
\end{center}}
\caption{Chosen parameters for the gamma example with $K=6$.}
\label{parameter values}
\end{table}

\begin{figure}[htb!]
\begin{center}
\begin{minipage}[t]{0.32\textwidth}
\begin{center}
\includegraphics[width=\textwidth]{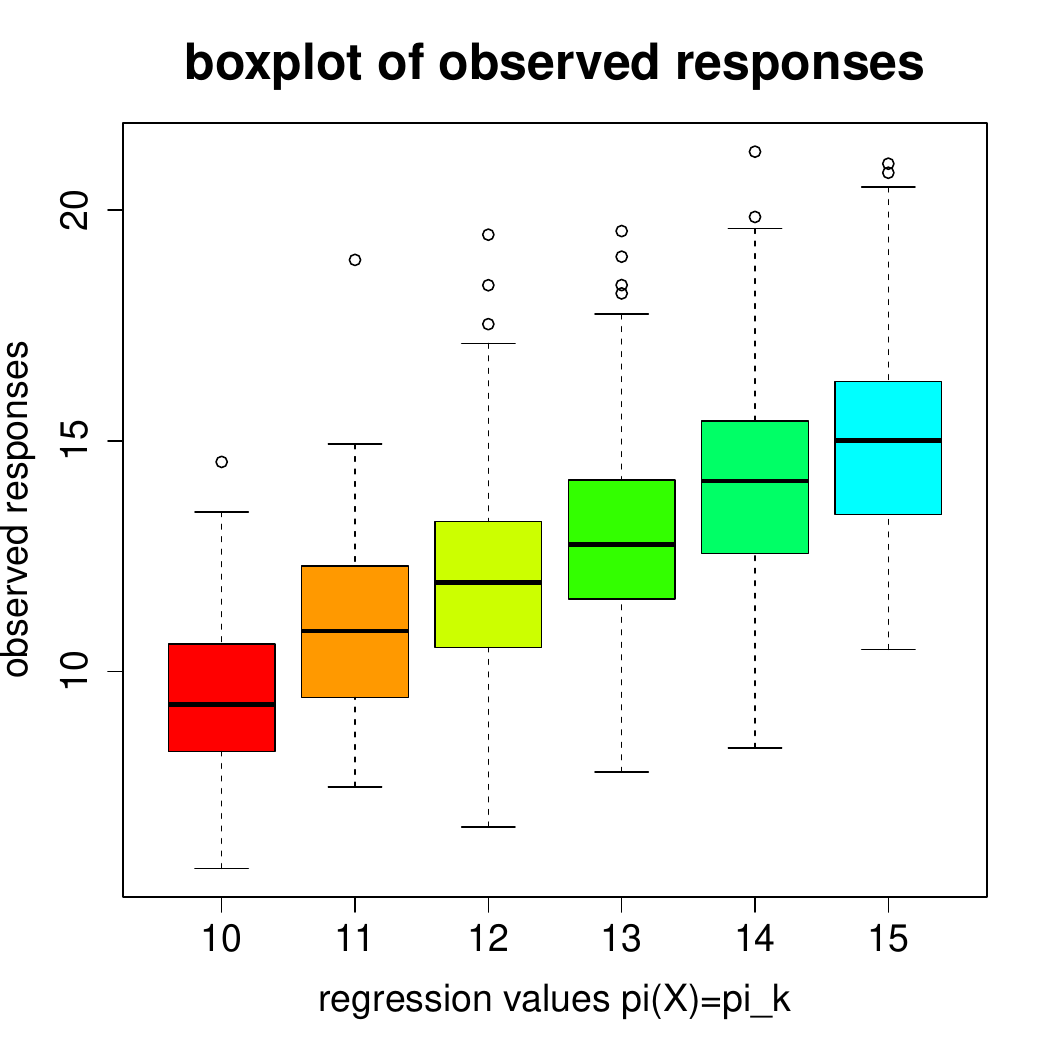}
\end{center}
\end{minipage}
\begin{minipage}[t]{0.32\textwidth}
\begin{center}
\includegraphics[width=\textwidth]{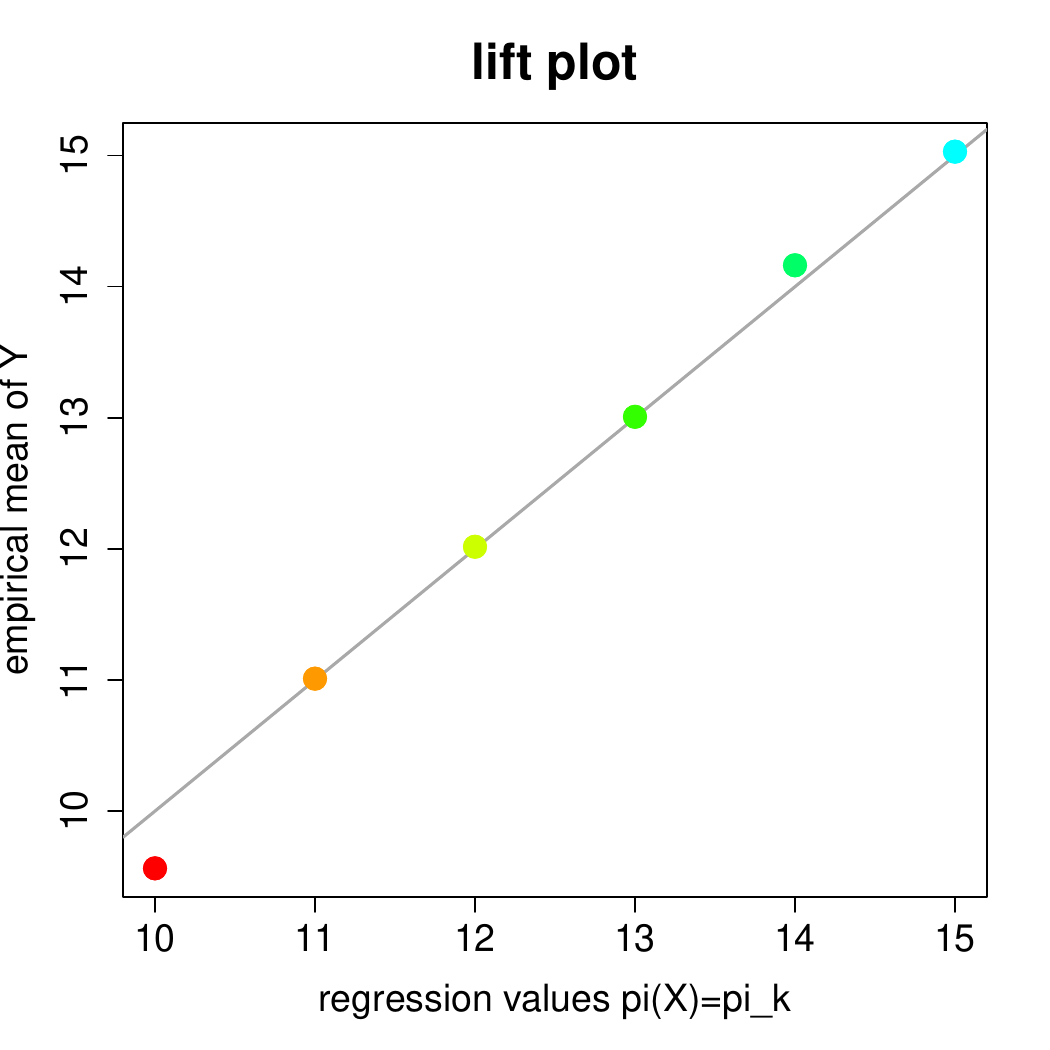}
\end{center}
\end{minipage}
\begin{minipage}[t]{0.32\textwidth}
\begin{center}
\includegraphics[width=\textwidth]{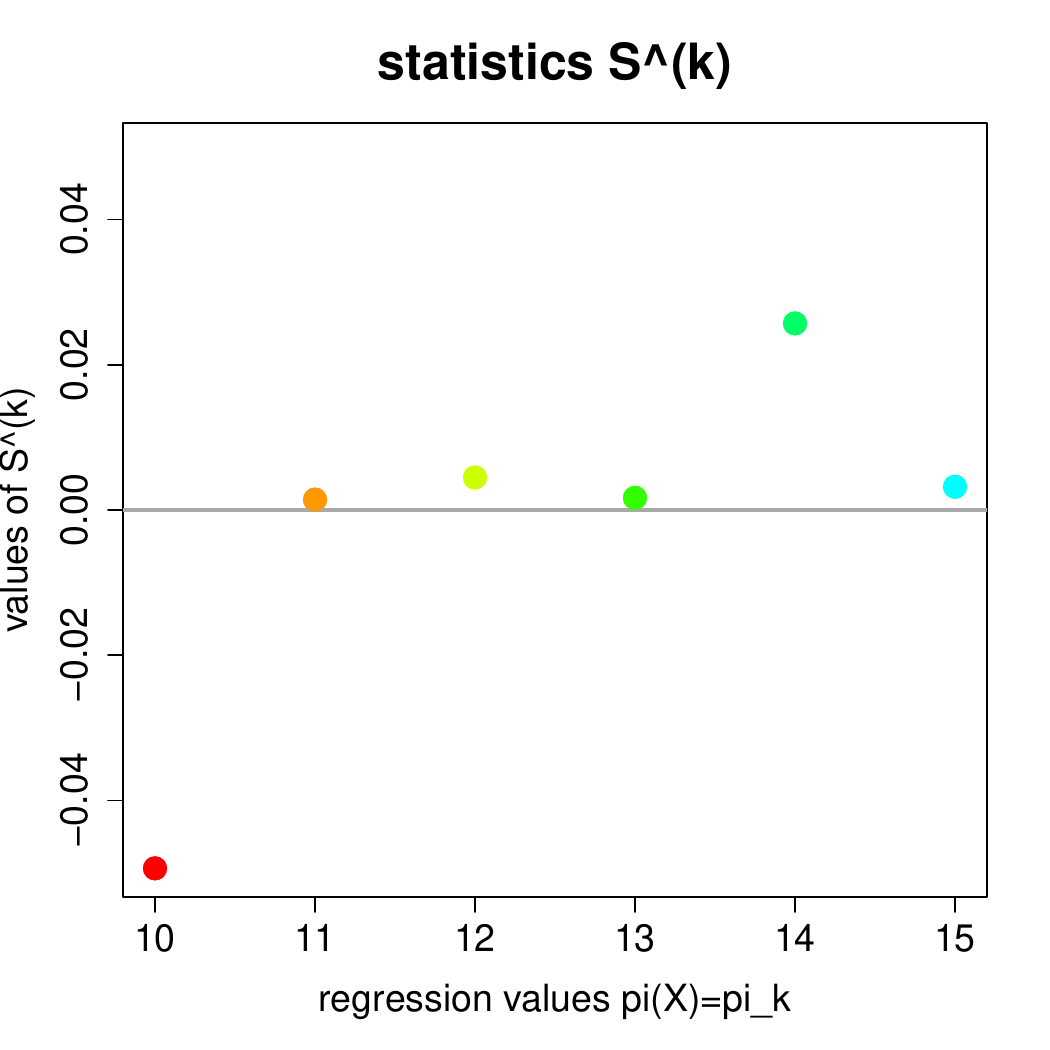}
\end{center}
\end{minipage}
\end{center}
\vspace{-.7cm}
\caption{Simulation of an i.i.d.~sample $(Y_i,\pi(\bX_i))_{i=1}^n$ of sample size $n=1000$: (lhs) boxplot of the responses $(Y_i)_{i=1}^n$ classified w.r.t.~$\pi(\bX_i)=\pi_k$, (middle) lift plot showing the
empirical level means $\overline{y}_k$ against their
expectations $\pi_k$, and (rhs) statistics $S^{(k)}_n$ for $1\le k \le K$.}
\label{Figure 1}
\end{figure}

Based on the parameters given in Table \ref{parameter values},
we simulate an i.i.d.~sample $(Y_i,\pi(\bX_i))_{i=1}^n$ of sample size $n=1000$. Figure \ref{Figure 1} (lhs) shows the resulting boxplot of the responses $(Y_i)_{i=1}^n$ classified w.r.t.~their conditional means $\pi(\bX_i)=\pi_k$. Remark that there is auto-calibration in this example.
Figure \ref{Figure 1} (middle) plots the empirical level means
\begin{equation*}
\overline{y}_k = \frac{1}{\sum_{i=1}^n\mathds{1}_{\{\pi(\bX_i)=\pi_k\}}}
\sum_{i=1}^n Y_i \,\mathds{1}_{\{\pi(\bX_i)=\pi_k\}},
\end{equation*}
against their (true) conditional expectations $\pi_k$; this plot is 
sometimes also called lift plot.
Under auto-calibration, the resulting scatter plot should lie fairly much on the
diagonal, and their deviation from the diagonal is described (asymptotically) by Proposition \ref{lemma 1}. Figure \ref{Figure 1} (rhs)
shows the resulting statistics $S^{(k)}_n$, for $1\le k \le K$.
These are obtained from the lift plot by using a different normalization
\begin{equation*}
S^{(k)}_n=\frac{1}{n} \sum_{i=1}^n \left(Y_i-\pi(\bX_i) \right) \mathds{1}_{\{\pi(\bX_i)=\pi_k \}} = \frac{\sum_{i=1}^k \mathds{1}_{\{\pi(\bX_i)=\pi_k\}}}{n}\left(\overline{y}_k - \pi_k\right),
\end{equation*}
the ratio on the right-hand side is an empirical estimate of $p_k$.
The magnitude of fluctuations of these statistics $S^{(k)}_n$ around
zero should be of order $\sqrt{p_k}\tau_k/\sqrt{n}$, see
Proposition \ref{lemma 1}.

\medskip

We repeat this simulation of an i.i.d.~sample $(Y_i,\pi(\bX_i))_{i=1}^n$ $10,000$ times to study the empirical distribution of the statistics
$\sqrt{n}\bS_n=\sqrt{n}(S^{(1)}_n, \ldots, S^{(K)}_n)^\top$. 
For large sample sizes $n$, this 
empirical distribution should approximately look
like the Gaussian limiting distribution
given in  Proposition \ref{lemma 1}.
Our simulation has an empirical mean $\widehat{\E}[\sqrt{n}\bS_n]$
of magnitude $10^{-2}$, thus, close to zero. 
The empirical covariance matrix reads as
\begin{equation*}
\widehat{\rm Cov}( \sqrt{n}\bS_n)=
\begin{pmatrix*}[r]
\bl{0.34} & 0.00 & 0.00 & -0.01 & 0.00 & 0.00\\
0.00 & \bl{0.55} &0.01 & 0.01 &0.00 &0.01\\
0.00& 0.01 & \bl{1.01} &-0.01 &0.00&  -0.01\\
-0.01 & 0.01 &-0.01 & \bl{1.10} & -0.01 & 0.01\\
0.00 &0.00& 0.00 & -0.01 & \bl{0.69} & 0.00\\
0.00& 0.01 & -0.01 & 0.01 & 0.00 & \bl{0.50}
\end{pmatrix*}.
\end{equation*}
The off-diagonals are close to zero and the diagonal is close
to true parameters
\begin{equation*}
\left(p_1 \tau_1^2, 
\ldots, p_5 \tau_5^2\right)
=\left( 0.33,\, 0.55,\, 1.00,\, 1.08,\, 0.70,\, 0.50\right),
\end{equation*}
see asymptotic covariance matrix in Proposition \ref{lemma 1}.
This confirms the limiting parameters in the weak convergence result of Proposition \ref{lemma 1}.

\begin{figure}[htb!]
\begin{center}
\begin{minipage}[t]{0.48\textwidth}
\begin{center}
\includegraphics[width=\textwidth]{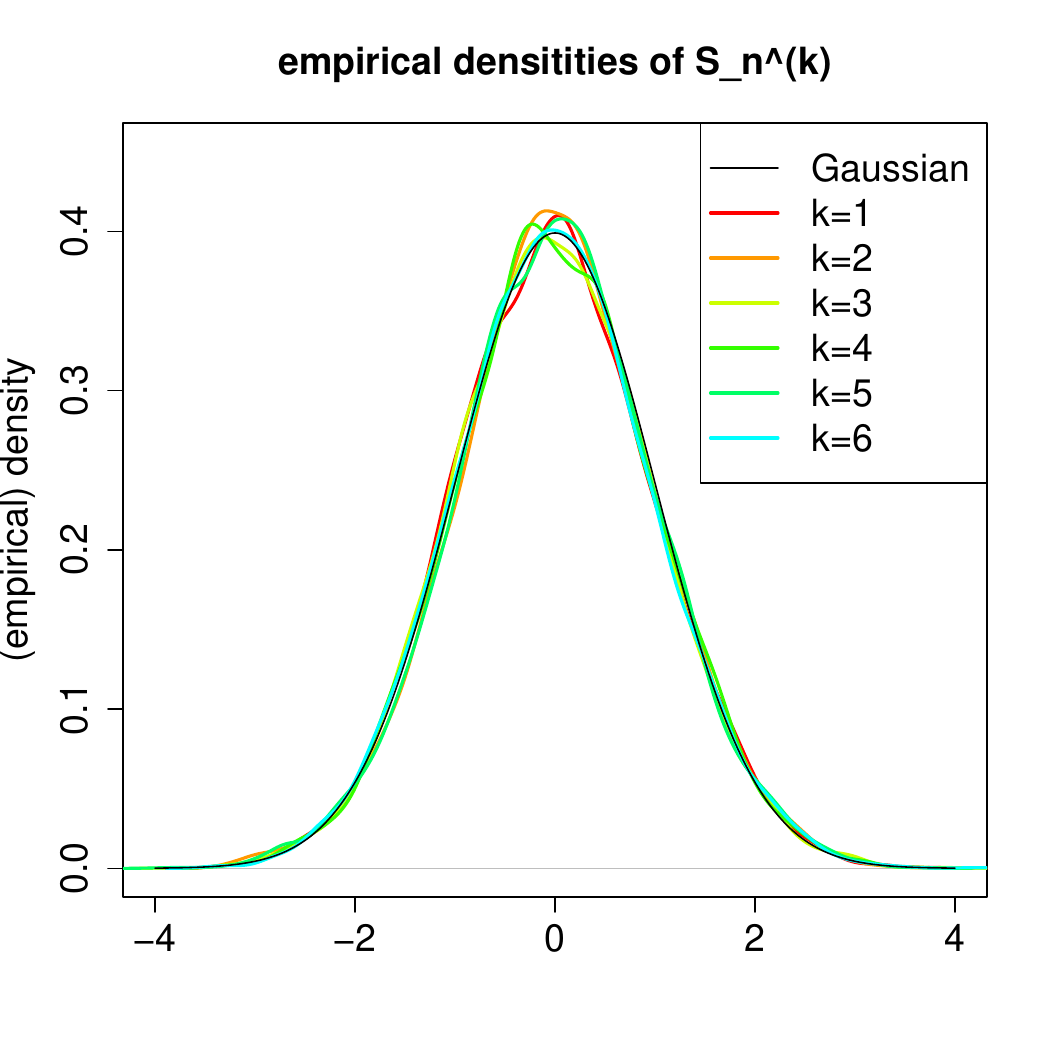}
\end{center}
\end{minipage}
\begin{minipage}[t]{0.48\textwidth}
\begin{center}
\includegraphics[width=\textwidth]{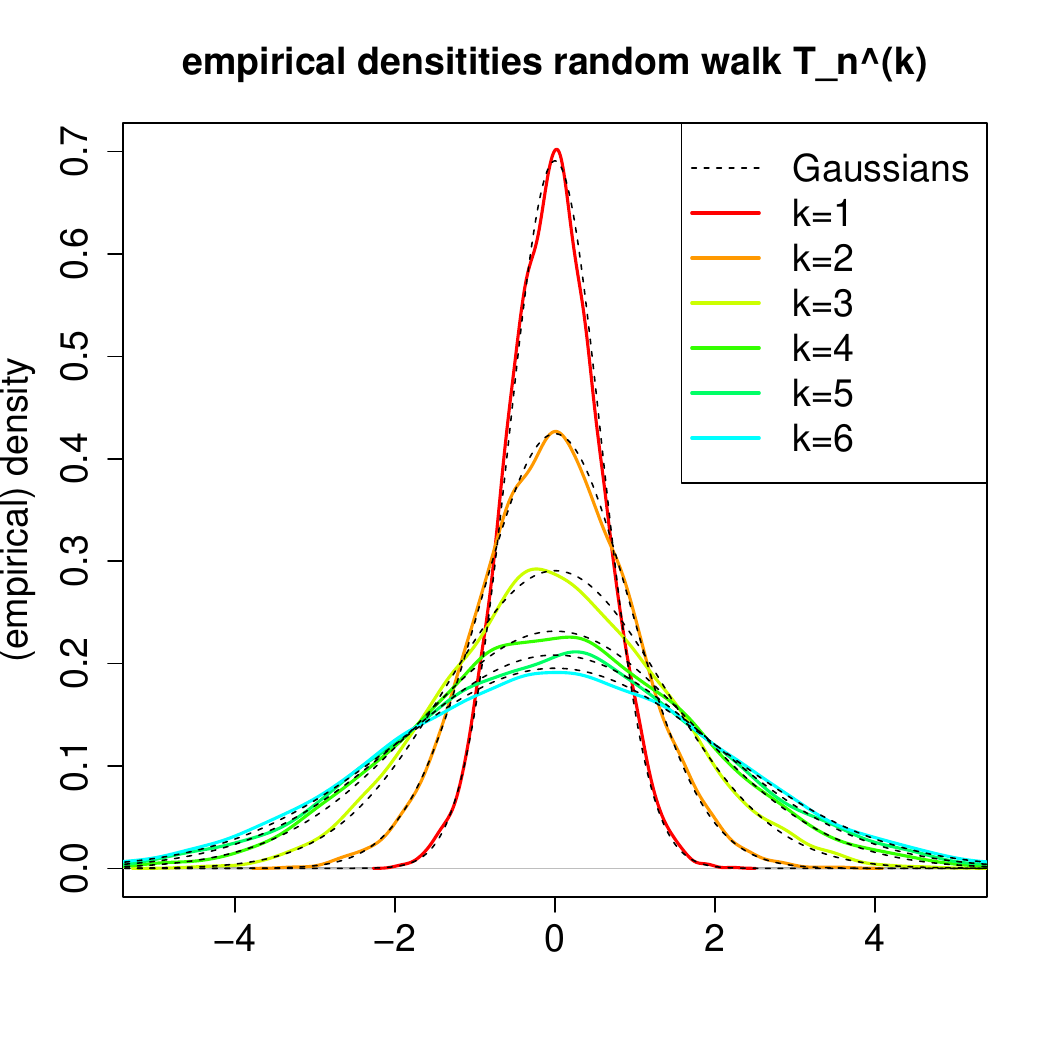}
\end{center}
\end{minipage}
\end{center}
\vspace{-.7cm}
\caption{(lhs) Empirical densities of $\sqrt{n}S^{(k)}_n/(\sqrt{p_k}\tau_k)$, for $1\le k \le K$, compared to the standard Gaussian density, 
and (rhs) empirical densities of the random walk $\sqrt{n}T^{(k)}_n$, for $1\le k \le K$, compared to the Gaussian random walk densities of $(Z_k)_{k=1}^K$, see \eqref{random walk}.}
\label{Figure 2}
\end{figure}

Figure \ref{Figure 2} (lhs) shows the empirical densities of $\sqrt{n}S^{(k)}_n/(\sqrt{p_k}\tau_k)$, for $1\le k \le K$ and sample size $n=1000$. They are benchmarked against the standard Gaussian density in black color. We see a quite good alignment of these empirical densities, supporting
the statement of Proposition \ref{lemma 1}. This justifies using the
asymptotic approximation \eqref{test 1} for the auto-calibration Test 1
in this example. Since the components in the maximum in
\eqref{test 1} may live on different scales, we also use an
alternative test statistics that evaluates the normalized quantities
\begin{equation}\label{test 1 B}
\p \left[ \max_{1\le k \le K}\sqrt{n} \left|\frac{S^{(k)}_n}{\sqrt{p_k}\tau_k}\right| \le s \right] = \p \left[ \bigcap_{1\le k \le K} 
\left\{\sqrt{n}\left|\frac{S^{(k)}_n}{\sqrt{p_k}\tau_k}\right| \le s\right\}
\right]  ~\approx~ \prod_{k=1}^K  \left(2 \Phi(s)-1\right).
\end{equation}
This then directly relates to the (normalized) graphs in Figure \ref{Figure 2} (lhs).

\medskip

Next, we turn our attention to the second test, involving the
random walk consideration \eqref{random walk}. In this case
we get the random walk type 
empirical covariance matrix
\begin{equation*}
\widehat{\rm Cov}( \sqrt{n}\bT_n)=
\begin{pmatrix*}[r]
\bl{0.34} &0.34 &0.34 &0.33 &0.33 &0.34\\
0.34 &\bl{0.89} &0.90 &0.89 &0.90 &0.91\\
0.34 &0.90 &\bl{1.91} &1.90 &1.91 &1.92\\
0.33 &0.89 &1.90 &\bl{2.99} &2.99 &3.01\\
0.33 &0.90 &1.91 &2.99 &\bl{3.68} &3.70\\
0.34 &0.91 &1.92 &3.01 &3.70 &\bl{4.21}
\end{pmatrix*}.
\end{equation*}
Since we work with a small sample size of $n=1000$, there is still some
noise involved which makes to above empirical covariance matrix not a perfect
random walk covariance matrix. The random walk covariance
matrix of Corollary \ref{corollary 2} has diagonal entries
$(0.33, 0.88, 1.88, 2.97, 3.67, 4.17)$.
Figure \ref{Figure 2} (rhs) plots the empirical densities $\sqrt{n}T^{(k)}_n$, for $1\le k \le K$, and these are benchmarked against the Gaussian random walk densities \eqref{random walk} of $Z_k$, $1\le k \le K$. Again we see a rather good alignment, supporting the asymptotic approximation \eqref{test 2} for auto-calibration Test 2. Clearly, the last random walk components
$\sqrt{n}T^{(K)}_n$ and $Z_K$, respectively, have the biggest variance, which implies that they will frequently determine the test statistics, see
\eqref{test 2}. Naturally, one could also revert index $k$
by studying the mirrored quantity, see also W\"uthrich (2023) for mirroring,
\begin{equation}\label{test 2 B}
\widetilde{T}^{(k)}=\E \left[ \left(Y-\pi(\bX) \right) \mathds{1}_{\{\pi(\bX)\ge \pi_k \}}\right]=0,
\end{equation}
and its empirical counterpart
\begin{equation*}
\widetilde{T}^{(k)}_n=\frac{1}{n} \sum_{i=1}^n \left(Y_i-\pi(\bX_i) \right) \mathds{1}_{\{\pi(\bX_i) \ge \pi_k \}}=
 \sum_{j=k}^K S^{(j)}_n.
\end{equation*}
If the terms $p_k \tau_k^2$ are increasing in $1\le k \le K$, this latter option may give a test with a better power, because the random walk increments will have a decreasing standard deviation.

\medskip

\begin{figure}[htb!]
\begin{center}
\begin{minipage}[t]{0.48\textwidth}
\begin{center}
\includegraphics[width=\textwidth]{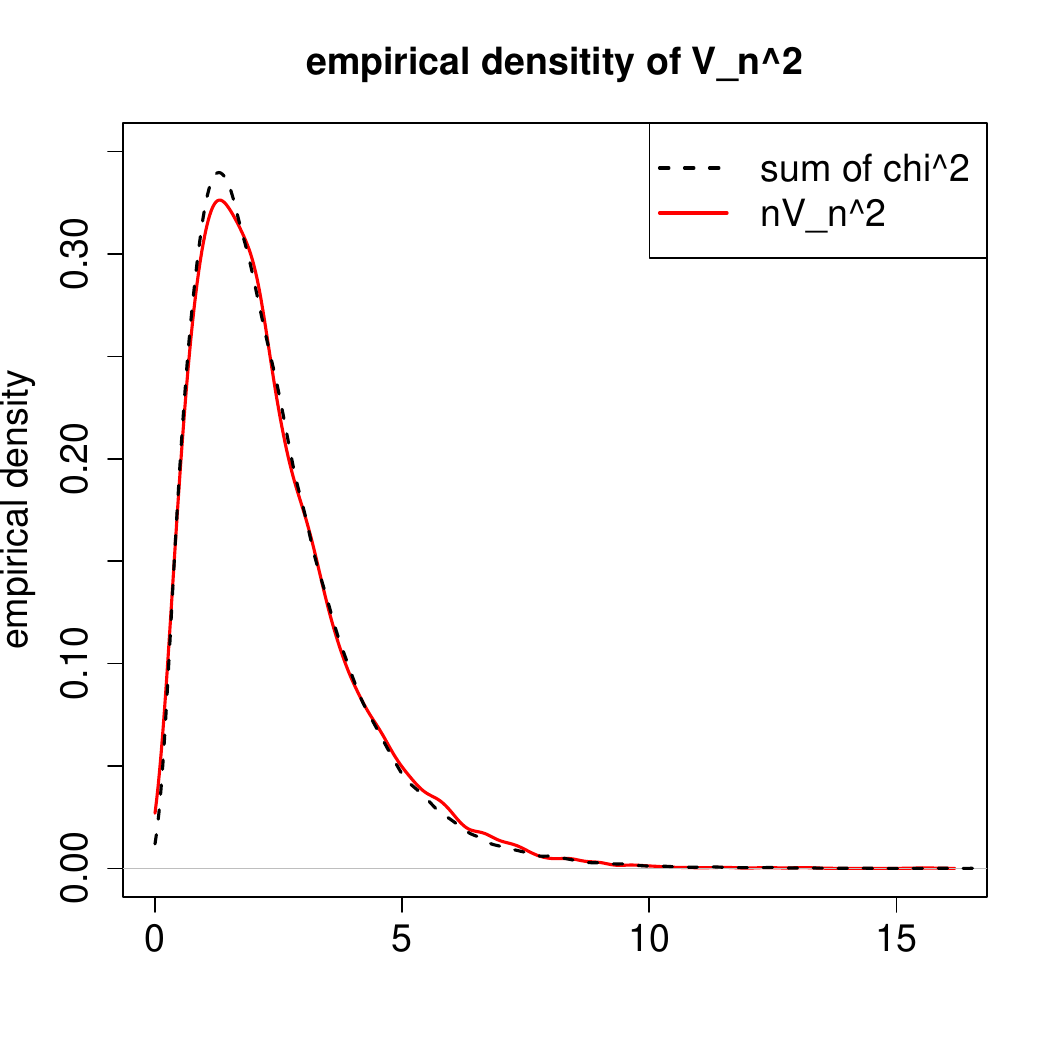}
\end{center}
\end{minipage}
\end{center}
\vspace{-.7cm}
\caption{Empirical density of $n V^2_n$ compared to the 
sum of independent scaled $\chi^2$-distributions given in 
Corollary \ref{corollary 3}.}
\label{Figure 3}
\end{figure}

Finally, Figure \ref{Figure 3} illustrates the asymptotic result
of Corollary \ref{corollary 3} for a sample size of $n=1000$. The test
statistics $nV_n^2$ does not consider the maximum of the increments,
$\max_{1\le k \le K}\sqrt{n}|S_n^{(k)}|$, but it considers a
weighted $L^2$-norm of all increments.
In \eqref{test 3} we study a weighted $L^2$-norm which has been
motivated by the ABC. However, in general, it is not clear why
this weighting should be justified. Alternatively, we could also
consider an unweighted test statistics
\begin{equation}\label{test 3B}
n\,\widetilde{V}^2_n
=n\, \sum_{k=1}^{K} (S_n^{(k)})^2\quad
 \Longrightarrow \quad \sum_{k=1}^{K} p_k\,\tau_k^2\,\chi_k^2
\qquad \text{ as $n\to \infty$},
\end{equation}
assuming $\pi$ is auto-calibrated for $(Y,\bX)$.
Equivalently, we could just consider a $\chi^2$-test
\begin{equation}\label{test 4}
n\, \sum_{k=1}^{K} \frac{(S_n^{(k)})^2}{p_k \, \tau^2_k}
\quad
 \Longrightarrow \quad \chi_K^2
\qquad \text{ as $n\to \infty$},
\end{equation}
where the right-hand side is a $\chi^2$-distributed random variable with $K$ degrees
of freedom. This is the same scaling as in \eqref{test 1 B}, however, we
do not consider maximums of increments, but rather aggregated squares of
the normalized random walk increments.

\bigskip

Summarizing, we have seen seven different test statistics that we will exploit numerically:
\begin{itemize}
\item[(1a)] From Test 1, we can study the maximum of the increments, see \eqref{test 1}.
\item[(1b)] A differently scaled version of Test 1 is given in \eqref{test 1 B}.
\item[(2a)] From Test 2, we can study the maximum of a random walk, see \eqref{test 2}.
\item[(2b)] An index reverted version of Test 2 is given in \eqref{test 2 B}.
\item[(3a)] From Test 3, we get a weighted $L^2$-norm of the
random walk increments, 
see \eqref{test 3}.
\item[(3b)] An unweighted alternative of Test 3 is
given in \eqref{test 3B}.
\item[(3c)] Finally, we have $\chi^2$-test given by \eqref{test 4}.
\end{itemize}
Because we have a discrete regression function $\pi$ taking finitely many
values, we receive a natural partition of the covariates space, ${\cal X}=\bigcup_{k=1}{\cal X}_k$,
and of the range of the regression function, $(\pi_k)_{k=1}^K$. For continuous
regression functions $\pi$, one can discretize the range of the regression function
$\pi$ and then perform a $\chi^2$-test for auto-calibration. In the Bernoulli case this has been proposed by Hosmer--Lemeshow (1980), and the
discretization is done with the help of the (empirical) quantiles of $\pi$. Our 
proposal is a generalization to arbitrary responses, and we present test
statistics that are different (and differently aggregated and normalized) from 
the classical $\chi^2$-test in the Bernoulli case.

\bigskip

Next, we aim at comparing the resulting powers of the seven tests in a simulation analysis.
We therefore contaminate the above model. We simulate
responses 
\begin{equation}\label{global shift}
Y^{\delta}= Y + \delta, \quad \text{ with }\quad
Y|_{\pi(\bX)=\pi_k} \sim \Gamma(\gamma_k,c_k).
\end{equation}
Thus, we introduce a global bias by shifting the means $\pi_k 
\mapsto \pi_k + \delta$ by a positive constant $\delta\ge 0$. 
This is a global shift as it affects equally all levels $\pi_k$, $1\le k \le K$.

\begin{table}[htb!]
\centering
{\small
\begin{center}
\begin{tabular}{|l||ccccccc|}
\hline
& Test 1a & Test 1b & Test 2a & Test 2b & Test 3a& Test 3b& Test 3c\\\hline
\hline
95\% quantiles  & 2.3456 & 2.6310 & 4.2060 & 4.2263 & 5.4066& 9.1198 & 12.5916\\
\hline
\end{tabular}
\end{center}}
\caption{Quantiles of the different test for significance level 5\%.}
\label{quantiles}
\end{table}

Table \ref{quantiles} gives the quantiles for significance level
5\% for the different tests.  The quantiles of Tests 1b and 3c are directly available in standard
software, the quantile of Test 1a can be found by a root search algorithm, and
quantiles of Tests 2a, 2b, 3a and 3b were computed empirically by a (simple) Monte Carlo simulation.

\begin{figure}[htb!]
\begin{center}
\begin{minipage}[t]{0.48\textwidth}
\begin{center}
\includegraphics[width=\textwidth]{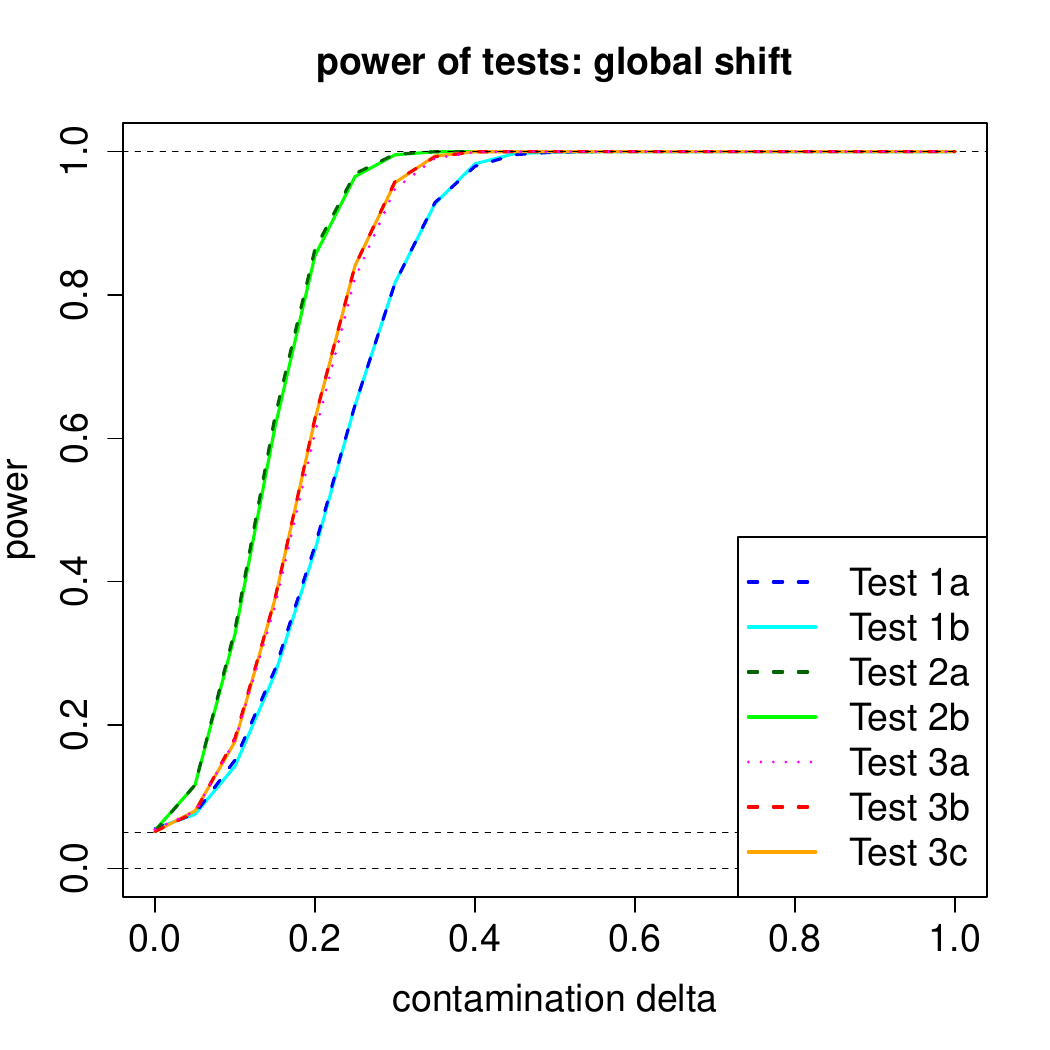}
\end{center}
\end{minipage}
\begin{minipage}[t]{0.48\textwidth}
\begin{center}
\includegraphics[width=\textwidth]{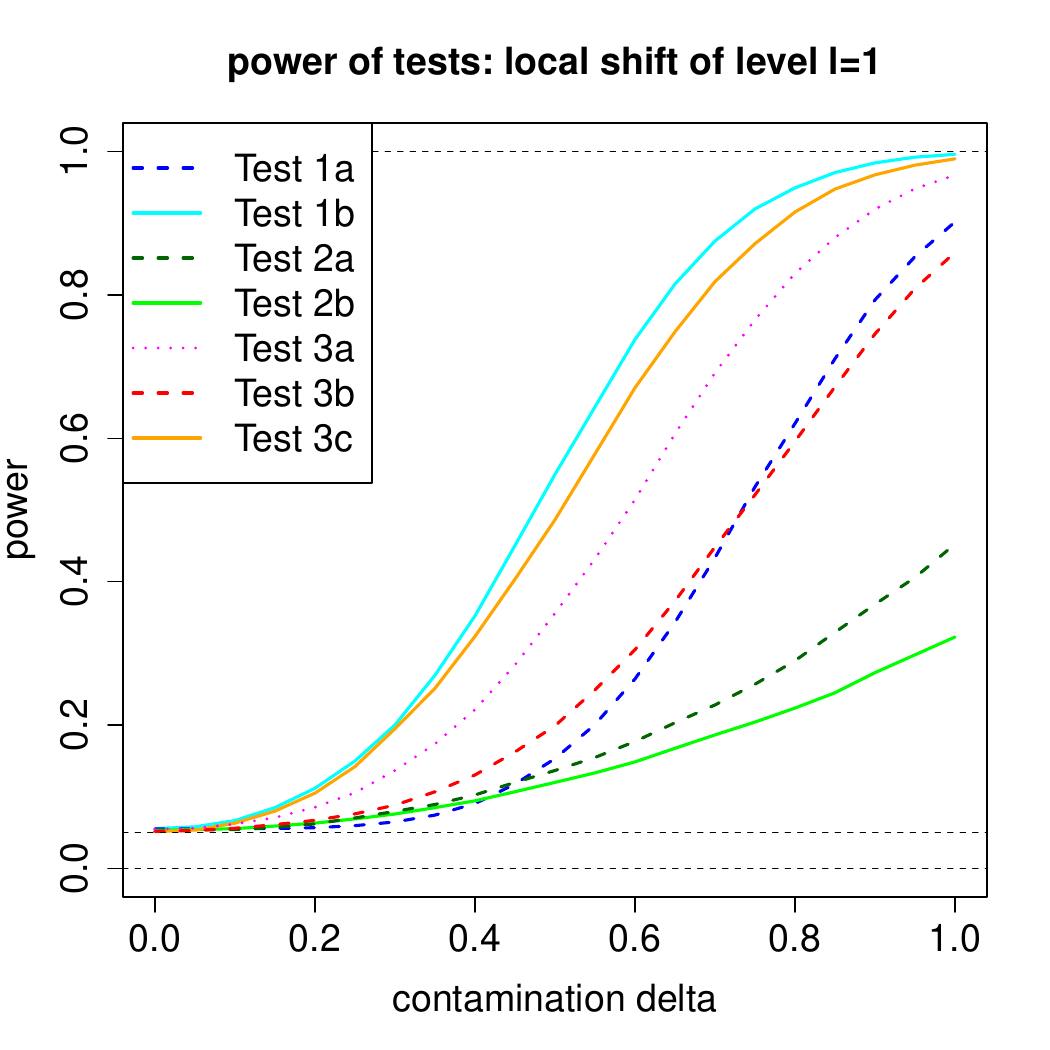}
\end{center}
\end{minipage}
\begin{minipage}[t]{0.48\textwidth}
\begin{center}
\includegraphics[width=\textwidth]{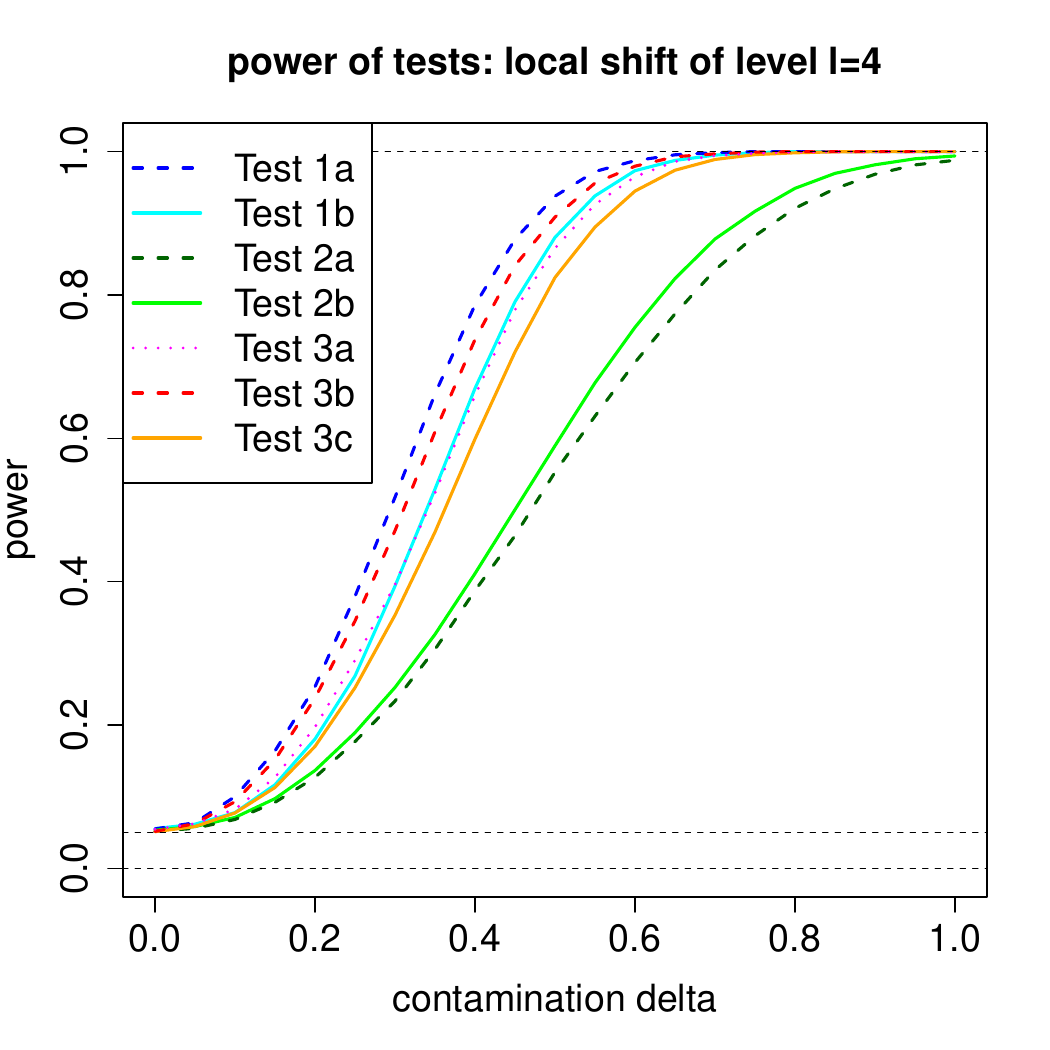}
\end{center}
\end{minipage}
\begin{minipage}[t]{0.48\textwidth}
\begin{center}
\includegraphics[width=\textwidth]{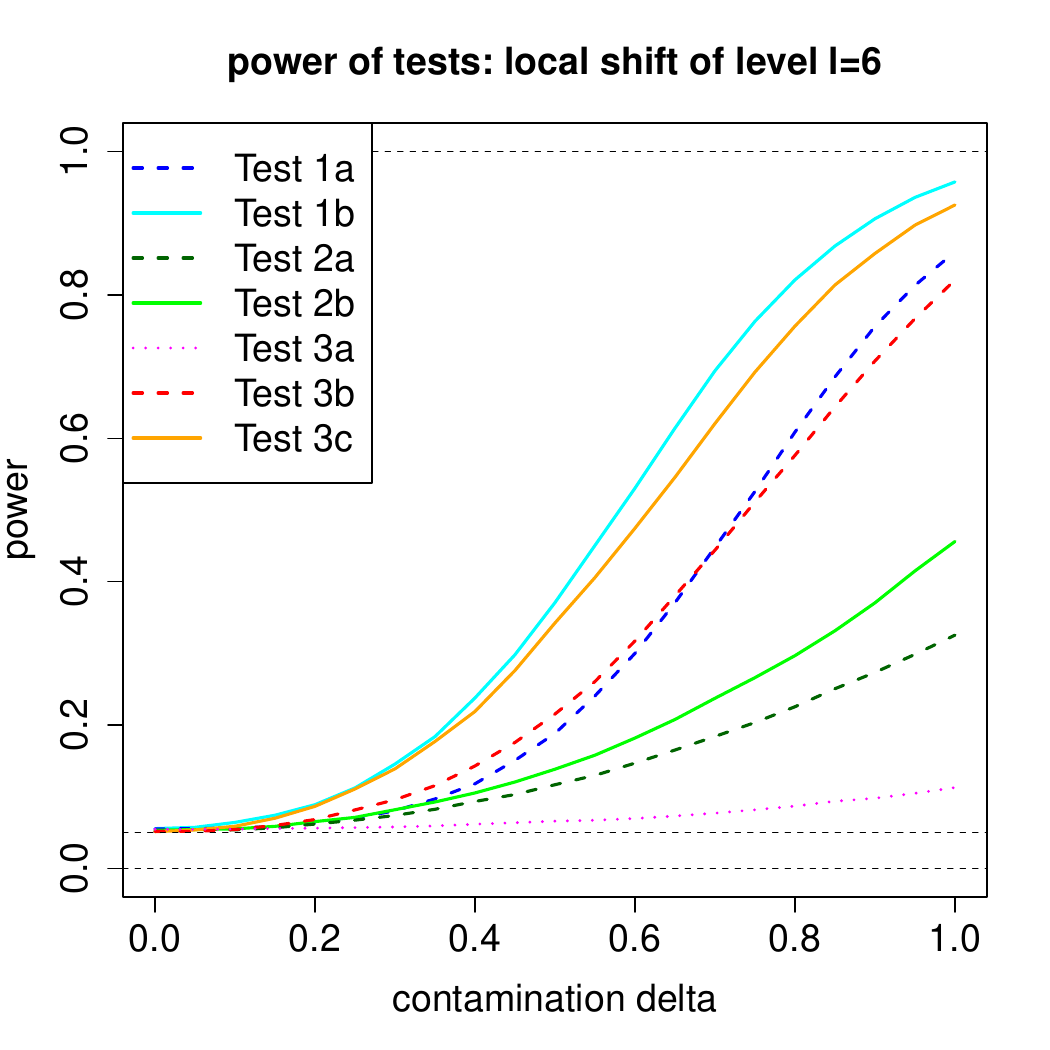}
\end{center}
\end{minipage}
\end{center}
\vspace{-.7cm}
\caption{Powers of the seven tests: (top-lhs) global contamination \eqref{global shift}, (top-rhs) local contamination \eqref{local shift} of 
the lowest level $\pi_1$, 
(bottom-lhs) local contamination \eqref{local shift} of level $\pi_4$,
(bottom-rhs) local contamination \eqref{local shift} of 
the highest level $\pi_6$.}
\label{Figure 4}
\end{figure}

We simulate 10,000 times (with different seeds)  i.i.d.~samples $(Y^\delta_i, \bX_i)_{i=1}^n$, $n=1000$, and for a grid of contaminations $\delta \in \{0, 1/20, 2/20, \ldots, 1\}$, see \eqref{global shift}.
This gives us for every simulation $1\le t \le 10,000$ and for every contamination level $\delta \in \{0, 1/20, 2/20, \ldots, 1\}$ the seven test
statistics. In the uncontaminated case $\delta=0$ roughly 5\% of the
10,000 simulations should be above the quantiles of Table \ref{quantiles}.
This then verifies that the asymptotic results for the tests apply, i.e., that $n=1000$
is a sufficiently large sample size for these tests.

For contaminations $\delta >0$ significantly more simulations should
be above the quantiles of Table \ref{quantiles}, and the more samples there are
above the corresponding quantile the bigger the power of the test. 
Figure \ref{Figure 4} (top-lhs) shows the results. We see that all curves
start at the significance level of 5\% for $\delta=0$. Then, they
increase to 1 for increasing contamination $\delta \uparrow 1$. The fastest increase
is achieved by Tests 2a-2b (maximum of random walk), followed by Tests 3a-3c (squared sum of random walk increments), and the slowest increase is achieved by Tests 1a-1b (maximum of random walk increments). From this we conclude that the
random walk tests \eqref{test 2} and \eqref{test 2 B} have the
biggest power in case of a global shift, and they should be preferred
to find global shifts. Intuitively this is clear, each random
walk increment $S_n^{(k)}$ is shifted by the contamination $\delta>0$,
and in the random walk these shifts are aggregated 
across all increments. Thus, we have an impact of $K\delta$ on the
last random walk component $T_n^{(K)}$. This is why Tests 2a-2b are the most sensitive ones to global shifts. In our example the order of aggregation is not very relevant, 
and Tests 2a-2b have almost equal power.

Global shifts are one potential cause of a violation of auto-calibration, but
the violation can also only occur on individual levels $\pi_k$, or on different
levels with different signs. To test for this local failure of the auto-calibration
property, we only contaminate individual levels of the regression function.
Fix $\ell \in \{1,\ldots, K\}$, and consider the
local contamination
\begin{equation}\label{local shift}
Y^{\delta, \ell}= Y + \delta\, \mathds{1}_{\{\pi(\bX)=\pi_\ell\}}, \quad \text{ with }\quad
Y|_{\pi(\bX)=\pi_k} \sim \Gamma(\gamma_k,c_k),
\end{equation}
this only contaminates the responses that have conditional
expectation $\pi(\bX)=\pi_\ell$. 

Based on this local contamination we repeat the above simulation experiment.
Since violation of auto-calibration often happens at the boundary of the range of the regression function, we contaminate
the model for the smallest and biggest conditional expectations $\pi_\ell$, $\ell \in \{1,6\}$. 
These are also the least frequent levels in our example. Additionally we contaminate level
$\pi_\ell$, $\ell=4$, being in the main body of the covariate distribution.
The results are presented in Figure \ref{Figure 4} (top-rhs and bottom).
The picture now significantly changes compared to the global contamination.
Tests 1b and 3c have the best behavior, both of these tests consider the normalized increments $S_n^{(k)}/(\sqrt{p_k}\tau_k)$. From this we conclude that one should bring
all random walk increments first to the same scale. This is especially true if the violation of auto-calibration takes place at  rare boundary levels, $\pi_1$ and $\pi_6$ in our case. For contaminated middle levels, $\pi_4$ in our case, the Tests 1a-1b and 3a-3c are all almost equally good. On the other hand, one should not use the aggregated random walk versions of Tests 2a-2b, because through aggregation the impact of individual violations of auto-calibration gets diluted. 
Another observation is that if the violation of auto-calibration happens on the biggest
level $\pi_6$, it cannot be found by the ABC inspired test
 \eqref{test 3}. This comes from the scaling $1-\alpha_{K-1}=p_K$ which
 often is a small number. Therefore, we cannot generally recommend Test 3a. 
 
 \medskip
 
 We summarize our findings of the simulation example as follows:
 \begin{itemize}
 \item Global shifts can most effectively be found by the random walk Tests 2a-2b, but this requires that auto-calibration is violated in the same direction
 on the entire support of the regression function.
 \item Local violation of auto-calibration, especially in the tails of the regression function can most effectively be found by Tests 1b and 3c. Both tests 
 consider scaled random walk increments (with unit variance), i.e., it seems beneficial that all random walk increments live on the same scale.
 \item The ABC inspired Test 3a can generally not be recommended, because the ABC weighting seems to prefer the lower over the upper tail of the regression function, but there
 is no specific reason that justifies such a weighting, compare magenta dotted lines in Figures \ref{Figure 4} (top-rhs) and (bottom-rhs).
 \end{itemize}

\subsection*{References}
\begin{itemize}
\item[]
Hosmer, D.W., Lemeshow, S. (1980). Goodness of fit tests for the multiple logistic regression model. {\it Communications
in Statistics - Theory and Methods} {\bf 9}, 1043-1069.
\item[]
W\"uthrich, M.V. (2023). Model selection with Gini indices under auto-calibration. {\it European Actuarial Journal} {\bf 13/1}, 469-477.
\end{itemize}

\end{document}